\documentclass[12pt,reqno]{amsart}

\usepackage{amssymb}

\usepackage{tikz-cd}

\usepackage{amsxtra}
\usepackage{verbatim}
\usepackage{enumerate}
\usepackage{enumitem}
\usepackage[english]{babel}
\usepackage{bm}

\usepackage{setspace}
\usepackage[x11names]{xcolor}

\usepackage[T1]{fontenc}

\usepackage{graphicx}

\usepackage{geometry}
\usepackage{amscd,latexsym,amsthm,amsmath,amsxtra}
 \usepackage[colorlinks, linkcolor=DodgerBlue3, citecolor=DeepPink1]{hyperref}
\usepackage{enumerate}
\usepackage[all]{xy}
\providecommand{\U}[1]{\protect\rule{.1in}{.1in}}

\RequirePackage{amsmath}
\RequirePackage{amssymb}

\usepackage[OT2,T1]{fontenc}
\DeclareSymbolFont{cyrletters}{OT2}{wncyr}{m}{n}
\DeclareMathSymbol{\sha}{\mathalpha}{cyrletters}{"58}

\input cyracc.def

\usepackage{comment}

 \newtheorem{thm}{Theorem}[section]

 \newtheorem{lem}[thm]{Lemma}

 \theoremstyle{definition}
 
 \theoremstyle{remark}
 
 \theoremstyle{remark}
 \newtheorem{rem}[thm]{Remark}

\newtheorem*{theorem*}{Theorem}
\newtheorem*{proposition*}{Proposition}
\newtheorem*{lemma*}{Lemma}
\newtheorem*{corollary*}{Corollary}
\newtheorem*{question*}{Question}
\newtheorem*{conjecture*}{Conjecture}
\newtheorem*{claim*}{Claim}

\newtheorem*{introtheorem*}{Theorem}
\newtheorem*{introproposition*}{Proposition}
\newtheorem*{introlemma*}{Lemma}
\newtheorem*{introcorollary*}{Corollary}

\numberwithin{equation}{section}

 \newcommand{\N}{\textup{N}}

 \renewcommand{\ker}{\textup{Ker}}

 \newcommand{\Br}{\textup{Br}}

 \newcommand{\Gal}{\textup{Gal}}

 \renewcommand{\P}{\mathbb{P}}
 \newcommand{\Ind}{\textup{Ind}}

 \newcommand{\Q}{\mathbb{Q}}

 \newcommand{\Z}{\mathbb{Z}}
 \newcommand{\G}{\mathbb{G}}

 \newcommand{\br}{\mathrm{Br}}
 
 \newcommand{\HH}{\mathrm{H}}

 \newcommand{\Res}{\mathrm{Res}}
 \newcommand{\Cor}{\mathrm{Cor}}
 \newcommand{\iden}{\mathrm{Id}}

\begin{document}

\title{Vertical unramified Brauer groups of Galois normic bundles }

\author{Yufan Liu}

\address{Yufan Liu \newline School of Mathematical Sciences, \newline  University of Science and Technology of China, \newline 96 Jinzhai Road, 230026 Hefei, China}

\email{liuyufan@mail.ustc.edu.cn}


\begin{abstract}
We compute the vertical unramified Brauer group of the Galois normic bundles, which are given  by $\N_{K/k}(\mathbf{z})=P(x)$.
Our main result gives combinatorial formulas for the vertical unramified Brauer groups in terms of the Galois group structure of $K/k$ and the irreducible factors of $P(x)$.  
\end{abstract}

\maketitle

\section{Introduction}

In the study of arithmetic geometry, computing the Brauer group of a given variety is a problem of considerable interest and complexity. The process is often intricate and far from trivial. Nevertheless, for rationally connected varieties, it is a well-established fact that their Brauer group (modulo the constant part) is always finite.

As a special class of rationally connected varieties, normic bundle is defined over a field $k$ by the equation
\begin{equation}\tag{$\star$}\label{normic eq}
\N_{K/k}(\mathbf{z})=P(x),
\end{equation}
where $K$ is a finite \'etale $k$-algebra, $\text{N}_{K/k}$ denotes the norm map, and $\mathbf{z}$ is a $K$-variable. 
For a smooth and proper model of \eqref{normic eq}, Colliot-Th\'el\`ene conjectured that the set of its rational points is dense in its Brauer--Manin set in \cite{CT03}. 
To compute the Brauer--Manin set of such a variety, a prerequisite is the computation of its unramified Brauer group.

The computation of the Brauer group for normic bundles has a rich history. In 2003, Colliot-Th\'el\`ene, Harari, and Skorobogatov constructed a partial compactification of the variety defined by \eqref{normic eq} and provided formulas for its vertical Brauer group and the corresponding quotient group in  \cite{CTHS03}.
While they successfully computed the unramified Brauer group for certain special cases (e.g., when $P(x)$ splits over $k$ with two distinct roots of relatively prime multiplicities), the Brauer group of their partial compactification is in general larger than the unramified Brauer group of the variety itself (see \cite[Theorem 3.6 and Proposition 3.7]{Wei14}). Building on the work of \cite{CTHS03}, Dasheng Wei investigated the case where $P(x)$ is irreducible and $K/k$ is abelian in \cite{Wei14}, while Derenthal, Smeets, and Wei studied the case where $\deg(P(x)) = 2$ and $[K : k] = 4$ in \cite{DSW12}. 

More recently, in \cite{wei2022+}, Wei constructed a partial compactification of normic bundles whose Brauer group coincides exactly with the unramified Brauer group. 
He characterized the vertical part of the Brauer group and its quotient via a short exact sequence of certain cohomology groups. 
Although these formulas hold for any finite \'etale $k$-algebra $K$ and any polynomial $P(x)$, they remain somewhat abstract in general.
For several special cases, Wei further determined the concrete group structures of the unramified Brauer group.

In parallel, in \cite{VAV-model}, V\'arilly-Alvarado and Viray succeeded in constructing explicit smooth and proper models of the normic bundle when $K/k$ is a cyclic extension of degree $n$, and $P(x)$ is a separable polynomial of degree $dn$ or $dn-1$.
In \cite{VAV-chatletp-fold}, they computed the Brauer group of such varieties when $n$ is prime.
Furthermore, in a recent paper \cite{LL25}, Yongqi Liang and the author computed the Brauer group of such normic bundles for general $n$ by using the same models. 
In the same paper, it was established that any non-zero finite abelian group can be realized as the Brauer group (modulo the constant part) of such a normic bundle.
However, constructing such smooth proper models in full generality remains a formidable task.

Building on the framework of \cite{wei2022+}, this paper provides the first systematic computation of the vertical part of the unramified Brauer group for a normic bundle corresponding to an arbitrary Galois extension $K/k$. 
As the main result of this paper, Theorem~\ref{main} explicitly characterizes this group by providing a concrete expression for its structure in all such cases.

The outline of the paper is as follows. In Section~\ref{notation}, we set up the notation and provide the preliminaries required for the subsequent cohomology calculations. In Section~\ref{review}, we review several results from \cite{wei2022+}, focusing in particular on the construction of the partial compactification. In Section~\ref{computation}, we carry out the explicit calculation of the vertical unramified Brauer group for Galois normic bundles.

\section{Notation and preliminaries}\label{notation}

\subsection{Notation and conventions}
In the present paper, we denote by $k$ a field of characteristic zero.  We fix an algebraic closure $\overline{k}$ of $k$. The absolute Galois group of $k$ is denoted by $\Gamma_k=\Gal(\overline{k}/k)$.

A variety over $k$ is understood as a separated scheme of finite type over $k$.

For any scheme $X$, its Brauer group is $\br(X):=\HH^2_\text{\'et}(X,\G_m)$.

For a field $K$ of characteristic zero, an element $a \in K^\times$, and a class $\chi \in \HH^1(K,\Q/\Z)$ whose order divides $m$, we write $a \smile \chi$ for the cup product of the image of $a$ in $K^\times / K^{\times m} = \HH^1(K,\mu_m)$ with $\chi$ in $\HH^1(K,\Z/m\Z)$, which lies in $\HH^2(K,\mu_m) = \Br(K)[m]\subset \br(K)$. 
In fact, the choice of $m$ does not matter here.

For a $\Gamma_k$-module $M$, we always denote
\[
    \sha_\omega^i(k,M):=\bigcap_{g\in \Gamma_k}\ker[\HH^i(k,M)\to  \HH^i(\langle g \rangle,M)].
\]

\subsection{preliminaries}
In this subsection, we will present some auxiliary results that will be used in the subsequent sections.

\begin{lem}[Shapiro's lemma]\label{shapiro}
    Let $G$ be a profinite group, $H$ an open subgroup of $G$, and $M$  a discrete $H$-module. 
    Then for any integer $i\geq 0$, we have
    \[
        \HH^i(G,\Ind_H^G M)\simeq \HH^i(H,M).
    \]

    Moreover, if $M$ is a discrete $G$-module with trivial $G$-action, then we can write $\Ind_H^G M$ as $\bigoplus g_i\otimes M$ where $g_i$ are representatives of the left cosets of $H$ in $G$.
    Define the map 
    \begin{align*}
      \phi:M&\to \Ind_H^G M,\\
      m&\mapsto \sum g_i\otimes m.
    \end{align*}
    Then for every $i\geq 0$, the composition of the following maps
    \[
       \HH^i(G,M)\xrightarrow{\phi_*} \HH^i(G,\Ind_H^G M)\xrightarrow{\simeq} \HH^i(H,M)
    \]
    coincides with the restriction map.
\end{lem}

\begin{lem}Let $G$ be a profinite group, $H$ an open subgroup of $G$, and $M$  a discrete $G$-module. 
  Assume that $G=\bigsqcup g_i H$. then the map
  \begin{align*}
      \Ind_H^G\Z\otimes M &\to \Ind^G_H(\Res^G_H M)\\
      (g_i\otimes 1)\otimes m&\mapsto g_i\otimes(g_i^{-1}m)
  \end{align*}
  is an isomorphism of $G$-modules.
\end{lem}

\begin{proof}
  This can be proved simply by verifying the $G$-action on both sides, since this map is clearly an isomorphism of abelian groups.
\end{proof}

\begin{lem}[Mackey's restriction formula]\label{mackey}
Let $G$ be a profinite group, $N$ and $H$ be two open subgroups, and $M$ a discrete $N$-module.
Then we have
\[
  \Res^G_H \Ind^G_N M \;\cong\;
  \bigoplus_{g \in H \backslash G / N}
  \Ind^{H}_{H \cap gNg^{-1}} \Res^{gNg^{-1}}_{H \cap gNg^{-1}} \, {}^g M,
\]
where ${}^g M$ denotes the $gNg^{-1}$-module obtained from $M$ by twisting the action via conjugation by $g$.

In particular, if $M$ is a discrete $G$-module with trivial $G$-action, and $N$ is a normal subgroup of $G$, then we have
\[
    \Res^G_H \Ind^G_N M \;\cong\;
  \bigoplus_{g \in H \backslash G / N}\Ind^{H}_{H \cap N} M.
\]  
\end{lem}

\begin{rem}
Throughout this paper, we make no distinction between induced modules and coinduced modules, as our considerations are restricted to open subgroups of profinite groups.
\end{rem}

\section{Review of Wei's results on normic bundles}\label{review}

To compute the unramified Brauer group of the normic bundle defined by the equation \eqref{normic eq}, Wei constructed a partial compactification of it in \cite{wei2022+}.
And he proved that the Brauer group of this partial compactification is just the unramified Brauer group of the original variety.
In this section, we will briefly review some of Wei's results in \cite{wei2022+}.

Let $k$ be a field of characteristic zero, $K/k$ be a Galois extension of degree $n$.
Let $P(x)\in k[x]$ be a polynomial of degree $m$.
Assume that $P(x)=c\cdot\prod_{i=1}^s P_i(x)^{e_i}  $, where $c\in k^\times$, $P_i(x)$ are monic irreducible polynomials in $k[x]$ and $e_i>0$ for each $i$. 
Let $U_1$ be the affine variety defined by the equation $\N_{K/k}(\mathbf{z})=P(x)$.

We will always assume that $n\mid m$, otherwise $U_1$ is birational to an affine variety of this form where the degree of the polynomial on the right side is divisible by $n$.
Similarly, we can always assume that $e_i<n$.

Let $U_2$ be the affine variety defined by the equation
\[
    \N_{K/k}(\mathbf{z}')=\widetilde{P}(x'),
\]
where $\widetilde{P}(x')={x'}^m\cdot P(1/x')$.
Then we can glue $U_1$ and $U_2$ along the open subset defined by $x\neq 0$ and $x'\neq 0$ respectively, by setting $x'=1/x$ and $\mathbf{z}'=\mathbf{z}/x^{m/n}$.
Let $U$ be the resulting variety.

Now we define $e_i'=\gcd(e_i,n)$ for $1\leq i\leq s$, and we may assume that $e_i'>1$ for $1\leq i\leq r$ and $e_i'=1$ for $r+1\leq i\leq s$.

Over $\overline{k}$, the affine variety $U_1$ can be defined by the equation:
\[
    P(x)=\prod_{j=1}^n z_j,
\]
where $z_i=\sigma_i \mathbf{z}$ for $\sigma_i\in \Gal(K/k)$.
Then for $1\leq i\leq r,1\leq j\leq n$, consider the closed subvariety of $\overline{U}$ (or $\overline{U_1}$) defined by
\[
    \bigcup_{1\leq i_1<i_2<\cdots<i_j\leq n}\left\{P_i(x)=z_{i_1}=z_{i_2}=\cdots z_{i_j}=0\right\}.
\]
One can easily check that this subvariety can be descended to a closed subvariety of $U$, we denote it by $W_{i,j}$.

Let $V:=U\setminus \bigcup_{i=1}^r W_{i,e_i'+1}$, and we blow $V$ along the ideal sheaf $I$ defined by the ideal:
\[
    \bigcap_{1\leq j\leq r} \ \bigcap_{1\leq i_1<i_2<\cdots<i_{e_j'}\leq n} (P_j(x)^{e_j/e_j'},z_{i_1},z_{i_2},\ldots,z_{i_{e_j'}}),
\]
where this ideal sheaf can be defined over $k$.
We denote the resulting variety by $\widetilde{V}$, and let $V_0$ be the smooth locus of $\widetilde{V}$.

Denote $U'=U\setminus \{P(x)=0\}$, then it is also an open subset of $V_0$.
One can check that any nonempty fiber of $\pi: U'\to \P^1$ is a torsor under the torus $T=\mathrm{R}_{K/k}^1(\G_m)$. 
Let $T^c$ be a fixed $T$-invariant smooth compactification of $T$, cf. \cite[Corollaire 1]{T^c}.
Then the contracted product $U'\times^T T^c$ contains $U'$ as an open subset.
We glue $V_0$ and $U'\times^T T^c$ along $U'$, and denote the resulting variety by $X$, this is just the partial compactification constructed by Wei in \cite{wei2022+}.
Then he proved the following theorem.

\begin{thm}[{\cite[Theorem 2.7]{wei2022+}}]
  The Brauer group of $X$ satisfies
  \[
      \br_1(X)=\br(X)=\br_{\mathrm{nr}}(X).
  \]
\end{thm}

\

Moreover, he also embedded $\br(X)/\br(k)$ into the middle term of an exact sequence.

Let $\widehat{T}$ be the character group of $T$, then $\widehat{T}\simeq \Z[K/k]/\Z$. 
Let $\Z_P:=\bigoplus \Z_{P_i}$ where each $\Z_{P_i}$ isomorphic to $\Ind_{\Gamma_{k(P_i)}}^{\Gamma_k}\Z$.
For $i=1,\ldots,r$, let $\Omega_i:=\{M:M\subset G,\#M=e_i' \}$, $S:=\bigoplus_{i=1}^r \Z_{P_i}\otimes \Z[\Omega_i]$.

  Then denote 
  \[
      \mathbb{D}:=\Z_P\otimes \Z[K/k]\oplus \Z\oplus S.
  \]

  Let $\N_K=\sum_{\sigma\in G}\sigma$ and $\Upsilon_i=\sum_{M\in \Omega_i}M$ for $i=1,\ldots,r$.
  We define a morphism $f:\Z_P\to \mathbb{D}$ by sending any $\tau\in \Gamma_k/\Gamma_{L_i}$ to $(\tau\otimes \N_K,-1,\tau\otimes \Upsilon_i)$ for $i=1,\ldots,r$ and sending any $\tau_i\in \Gamma_k/\Gamma_{L_i}$ to $(\tau\otimes \N_K,-1,0)$ for $i=r+1,\ldots,s$.
  Let $\widehat{T}'$ be the cokernel of $f$.

  Then we can define a map $j$ from $\widehat{T}$ to $\widehat{T}'$ which is induced by the map
  \begin{gather*}
    \Z[K/k]\to \mathbb{D},\\
    \sigma\mapsto (-\sum_{i=1}^s e_i \N_i\otimes \sigma,m/n,-S_\sigma)
  \end{gather*}
  where $\sigma\in G$, $\N_i=\sum_{\tau\in \Gamma_k/\Gamma_{L_i}}\tau$ and 
  \[
      S_\sigma=\sum_{i=1}^r\N_i \otimes (\frac{e_i}{e_i'}\cdot \sum_{M\in \Omega_i,\sigma\in M}M).
  \]
  
\begin{thm}[{\cite[Corollary 2.9 and Remark 3.17]{wei2022+}}]\label{exact}
  We have the exact sequence:
  \[
      0\to \HH^1(k,\widehat{T}')/j^*(\HH^1(k,\widehat{T}))\to \br(X)/\br(k) \to \ker[\sha_\omega^2(k,\widehat{T})\to\sha_\omega^2(\widehat{T}')]\to 0.
  \]

  Moreover, let $\br_{\mathrm{vert}}(X):=\br(X)\cap \br(k(x))$ be the vertical part of the Brauer group, then $\br_{\mathrm{vert}}(X)/\br(k)$ is just the first term of the above exact sequence, i.e. $\HH^1(k,\widehat{T}')/j^*(\HH^1(k,\widehat{T}))$.
\end{thm}

\begin{rem}
  In \cite{wei2022+}, Wei gave a similar exact sequence when $K/k$ is merely a finite \'etale algebra, rather than a Galois extension. Since our computations concern only the vertical Brauer group in the case where $K/k$ is a Galois extension, we do not review the general situation here, whose partial compactification is more complicated.
\end{rem}

\section{Computation of vertical Brauer groups}\label{computation}

In this section, we will compute the vertical part of the Brauer group of $X$ which is constructed in the previous section.

With the notation in the previous section, let $L_i/k$ be the field extension defined by $P_i(x)$ for $i=1,\ldots, s$, that is $L_i:=k[x]/(P_i(x))$.
Then $KL_i/L_i$ is still a  Galois extension for each $i$.  
Setting $L_i':=L_i\cap K$, then one can check that $\Gal(KL_i/L_i)=\Gal(K/L_i')$.
Let $l_i=[L_i:L_i']$.

Let $G:=\Gal(K/k)$, and for each $i$ set $G_i:=\Gal(KL_i/L_i)=\Gal(K/L_i')$. 
Then clearly $G_i\subset G$ for each $i$.
Now let $\widehat{G}:=\HH^1(G,\Q/\Z)$, $\widehat{G}_i:=\HH^1(G_i,\Q/\Z)$. 
Then for each $i$, there are natural maps $\Res_i:\widehat{G}\to\widehat{G}_i$ and  $\Cor_i:\widehat{G}_i\to\widehat{G}$, satisfying the relation $\Cor_i\circ\Res_i=[L_i':k]\cdot\iden_{\widehat{G}}$.

Now let $\widehat{G}_i':=\{\chi\in\widehat{G}_i:\chi(g)=0 \text{ for any } g\in G_i \text{ with } g^{e_i'}=1_{G_i}\}$, then we have the following theorem.

\begin{thm}\label{main}
  There is an isomorphism
    \[
        \br_{\textup{vert}}(X)/\br(k)\simeq \frac{\{(\chi_i)_{i=1}^s\in\bigoplus_{i=1}^s\widehat{G}_i':\sum_{i=1}^s l_i\cdot\Cor_i (\chi_i)=0\in \widehat{G}\}}{\{(e_i\cdot \Res_i(\chi))_{i=1}^s:\chi\in\widehat{G}\}}.
    \]
\end{thm}

\begin{proof}
  According to Theorem \ref{exact}, we only need to prove that 
  \[
      \HH^1(k,\widehat{T}')/j^*(\HH^1(k,\widehat{T}))\simeq \frac{\{(\chi_i)_{i=1}^s\in\bigoplus_{i=1}^s\widehat{G}_i':\sum_{i=1}^s l_i\cdot\Cor_i (\chi_i)=0\in \widehat{G}\}}{\{(e_i\cdot \Res_i(\chi))_{i=1}^s:\chi\in\widehat{G}\}}.
  \]

  By Lemma~\ref{shapiro}, one can easily check that $\HH^1(k,\Z_P)=\HH^1(k,\mathbb{D})=0$.
  Then $\HH^1(k,\widehat{T}')=\ker[\HH^2(k,\Z_P)\to \HH^2(k,\mathbb{D})]$.
  Since $\mathbb{D}=\Z_P\otimes \Z[K/k]\oplus \Z\oplus S$, we conclude that $\HH^1(k,\widehat{T}')$ is precisely the intersection of the kernels of the following maps 
  \begin{itemize}
    \item $\HH^2(k,\Z_P)\to \HH^2(k,\Z_P\otimes \Z[K/k])$,
    \item $\HH^2(k,\Z_P)\to \HH^2(k,\Z)$,
    \item $\HH^2(k,\Z_P)\to \HH^2(k,S)$.
  \end{itemize}

  Also by Lemma~\ref{shapiro}, we know that $\HH^2(k,\Z_P)=\bigoplus \HH^2(k,\Z_{P_i})=\bigoplus \HH^2(L_i,\Z)$.
  According to the long exact sequence from the short exact sequence
  \[
      0\to \Z \to \Q \to \Q/\Z \to 0,
  \]
  we have $\HH^2(L_i,\Z)=\HH^1(L_i,\Q/\Z)$ since  $\HH^j(L_i,\Q)=0$ for any $j>0$. 
  Then $\HH^2(k,\Z_P)=\bigoplus \HH^1(L_i,\Q/\Z)$.

  Similarly, we have 
  \[
      \HH^2(k,\Z_P\otimes \Z[K/k])=\bigoplus_{i=1}^s\HH^2(L_i,\Res_{\Gamma_{L_i}}^{\Gamma_k}\Z[K/k]))=\bigoplus_{i=1}^s\HH^2(L_i,\Res_{\Gamma_{L_i}}^{\Gamma_k}(\Ind_{\Gamma_K}^{\Gamma_k}\Z)).
  \]
  By Lemma~\ref{mackey}, we have $\Res_{\Gamma_{L_i}}^{\Gamma_k}(\Ind_{\Gamma_K}^{\Gamma_k}\Z)=\bigoplus \Ind_{\Gamma_{KL_i}}^{\Gamma_{L_i}}\Z$, where the number of summands of this direct sum is equal to the number of representatives of the double cosets $\Gamma_{K}\backslash \Gamma_k/\Gamma_{L_i}$.
  Then $\HH^2(k,\Z_P\otimes \Z[K/k])=\bigoplus_{i=1}^s\bigoplus \HH^1(KL_i,\Q/\Z)$, and the map from $\HH^1(L_i,\Q/\Z)$ to $\HH^1(KL_i,\Q/\Z)$ is just the restriction map by Lemma~\ref{shapiro}.
  Then we have 
  \[
      \ker[\HH^2(k,\Z_P)\to \HH^2(k,\Z_P\otimes \Z[K/k])]=\bigoplus \HH^1(\Gal(KL_i/L_i,\Q/\Z))=\bigoplus \widehat{G}_i.
  \]

  Secondly, we will compute $\ker [\HH^2(k,\Z_P)\to \HH^2(k,\Z)]$.
  By Lemma~\ref{shapiro}, we know that $\HH^2(k,\Z_P)=\bigoplus\HH^1(L_i,\Q/\Z)$ and $\HH^2(k,\Z)=\HH^1(k,\Q/\Z)$.
  Then one can check that the map from each summand $\HH^1(L_i,\Q/\Z)$ to $\HH^1(k,\Q/\Z)$ is just the corestriction map.
  Thus, we have the following commutative diagram
  \[
  \begin{tikzcd}
  \widehat{G}_i \arrow[d, hook] \arrow[rr, "\mathrm{Cor}_i"] &  & \widehat{G} \arrow[d, "\times l_i"] \\
  {\HH^1(L_i,\Q/\Z)} \arrow[rr, "\mathrm{Cor}_{L_i/k}"]      &  & {\HH^1(k,\Q/\Z),}                   
\end{tikzcd}
  \]
  which follows from the commutativity of the diagram below
  \[
  \begin{tikzcd}
\widehat{G}_i \arrow[d, hook] \arrow[r, "="]     & \widehat{G}_i \arrow[d, hook] \arrow[r, "="]    & \widehat{G}_i \arrow[d, "\times l_i"] \arrow[rr, "\Cor_i"] &  & \widehat{G} \arrow[d, "\times l_i"] \\
{\HH^1(L_i',\Q/\Z)} \arrow[r, "\Res_{L_i/L_i'}"] & {\HH^1(L_i,\Q/\Z)} \arrow[r, "\Cor_{L_i/L_i'}"] & {\HH^1(L_i',\Q/\Z)} \arrow[rr, "\Cor_{L_i'/k}"]            &  & {\HH^1(k,\Q/\Z).}                   
\end{tikzcd}
  \]
   It implies that 
  \[
      \ker [\HH^2(k,\Z_P)\to \HH^2(k,\Z)]\cap \bigoplus\widehat{G}_i=\{(\chi_i)_{i=1}^s\in \bigoplus\widehat{G}_i:\sum l_i\cdot \Cor(\chi_i)=0\}.
  \]

  Now we will prove that $\bigoplus \widehat{G}_i\cap \ker[\HH^2(k,\Z_P)\to \HH^2(k,S)]=\bigoplus \widehat{G}_i'$.
  By the definition of $f$, one can easily check that such group is just 
  \[
      \bigoplus_{i=1}^r (\widehat{G}_i\cap \ker[\HH^2(k,\Z_{P_i})\to \HH^2(k,\Z_{P_i}\otimes \Z[\Omega_i])])\oplus \bigoplus_{j=r+1}^s \widehat{G}_j.
  \]
  Since we have $\widehat{G}_j=\widehat{G}_j'$ for $j=r+1,\ldots,s$, we only need to prove that 
  \[
      \widehat{G}_i\cap \ker[\HH^2(k,\Z_{P_i})\to \HH^2(k,\Z_{P_i}\otimes \Z[\Omega_i])]=\widehat{G}_i'
  \]
  for $i=1,\ldots,r$. 
  By the argument above, we know that $\HH^2(k,\Z_{P_i})=\HH^1(L_i,\Q/\Z)$ and $\HH^2(k,\Z_{P_i}\otimes \Z[\Omega_i])=\HH^2(L_i,\Res_{\Gamma_{L_i}}^{\Gamma_k}\Z[\Omega_i])$.

  By the definition of $\Omega_i$, the group action of $\Gamma_k$ on $\Omega_i$ factor through $G$.
  Then the action of $\Gamma_{L_i}$ on $\Omega_i$ factor through $G_i$.
  In fact we have $\Omega_i=\bigsqcup_{j=1}^t G_i M_j$ for some $M_j\in \Omega_i$.
  Let $G_{i,j}$ be the stabilizer subgroup of $G_i$ respect to $M_j$, and $G_{i,M}$ be the stabilizer subgroup of $G_i$ respect to any $M\in \Omega_i$.
  Now let $L_{i,j}$ (resp. $L_{i,M}$) be the fixed field of $G_{i,j}$ (resp. $G_{i,M}$).
  As a $\Gamma_{L_i}$-module, we have the identifications $\Z[G_i M_j]\simeq\Ind_{\Gamma_{L_{i,j}}}^{\Gamma_{L_i}}\Z$ and $\Z[G_i M]\simeq\Ind_{\Gamma_{L_{i,M}}}^{\Gamma_{L_i}}\Z$.
  Hence, it follows that $\Res_{\Gamma_{L_i}}^{\Gamma_k}\Z[\Omega_i]=\bigoplus_{j=1}^t \Ind_{\Gamma_{L_{i,j}}}^{\Gamma_{L_i}}\Z$.

  According to Lemma~\ref{shapiro}, $\HH^2(k,\Z_{P_i}\otimes \Z[\Omega_i])=\bigoplus_{j=1}^t \HH^1(L_{i,j},\Q/\Z)$, and the map from $\HH^1(L_i,\Q/\Z)$ to each summand of this direct sum induced by $f$ is just the restriction map.
  Then for any $\chi\in \widehat{G}_i$, we have
  \begin{align*}
    &\chi\in \widehat{G}_i\cap \ker[\HH^2(k,\Z_{P_i})\to \HH^2(k,\Z_{P_i}\otimes \Z[\Omega_i])]\\
    \iff  & \chi|_{G_{i,j}}=0\textup{ for any }j=1,\ldots,t\\
    \iff & \chi|_{G_{i,M}}=0\textup{ for any } M\in \Omega_i.
  \end{align*}

  Since any $M\in \Omega_i$ is a subset of $G$ which has $e_i'$ elements, one can check that the exponent of $G_{i,M}$ is dividing $e_i'$ for any $M\in \Omega_i$.
  Also one can easily check that for any $g\in G_i$ of order dividing $e_i'$, there exists $M\in \Omega_i$ such that $g\in G_{i,M}$.
  Thus, based on the above equivalence relations, we have $ \widehat{G}_i\cap \ker[\HH^2(k,\Z_{P_i})\to \HH^2(k,\Z_{P_i}\otimes \Z[\Omega_i])]=\widehat{G}_i'$.

  By the above arguments, we have established that 
  \[
      \HH^1(k,\widehat{T}')=\{(\chi_i)_{i=1}^s\in\bigoplus_{i=1}^s\widehat{G}_i':\sum_{i=1}^s l_i\cdot\Cor_i (\chi_i)=0\in \widehat{G}\}.
  \]
  Moreover, by the definition of $j$, it follows that 
  \[
      j^*(\HH^1(k,\widehat{T}))=\{(e_i\cdot \Res_i(\chi))_{i=1}^s:\chi\in\widehat{G}\}.
  \]
  This completes the proof.
\end{proof}
\
\begin{rem}
    According to Lemma 2.5 in \cite{wei2022+}, we have an explicit description of the correspondence between $\bigoplus \widehat{G}_i'$ and $\br(k(x))$.
    Let $\varepsilon_i$ be the class of $x$ in $L_i=k[x]/(P_i(x))$, in fact it is a root of $P_i(x)$ in $L_i$.
    Then $(\chi_i)_{i=1}^s$ corresponds to $\sum_{i=1}^{s}\Cor_{L_i(x)/k(x)}((x-\varepsilon_i)\smile\chi_i)$ in $\br(k(x))$.
    In particular, if $\chi_i=\Res_i(\chi)$ for some $\chi\in \widehat{G}$, then by Proposition 1.5.3 in \cite{cnf}, we have 
    \[\Cor_{L_i(x)/k(x)}((x-\varepsilon_i)\smile\chi_i)=P_i(x)\smile \chi.\]
\end{rem}

\vspace{1cm}

{\bf Acknowledgements}
{\it The author would like to thank Professor Dasheng Wei for his valuable discussions and suggestions on this topic.}

\bibliographystyle{amsalpha}
\bibliography{mybib1}

\providecommand{\bysame}{\leavevmode\hbox to3em{\hrulefill}\thinspace}
\providecommand{\MR}{\relax\ifhmode\unskip\space\fi MR }
\providecommand{\MRhref}[2]{%
  \href{http://www.ams.org/mathscinet-getitem?mr=#1}{#2}
}
\providecommand{\href}[2]{#2}
\begin{thebibliography}{CTHS05}

\bibitem[CT03]{CT03}
J.-L. Colliot-Th{\'e}l{\`e}ne, \emph{Points rationnels sur les fibrations}, Higher Dimensional Varieties and Rational Points (K{\'a}roly B{\"o}r{\"o}czky, J{\'a}nos Koll{\'a}r, and Tam{\'a}s Szamuely, eds.), Springer Berlin Heidelberg, Berlin, Heidelberg, 2003, pp.~171--221.

\bibitem[CTHS04]{CTHS03}
J.-L. Colliot-Th{\'e}l{\`e}ne, D.~Harari, and A.~N. Skorobogatov, \emph{Valeurs d'un polyn{\^o}me {\`a} une variable repr{\'e}sent{\'e}es par une norme}, Number Theory and Algebraic Geometry (Miles Reid and Alexei Skorobogatov, eds.), London Mathematical Society Lecture Note Series, vol. 303, Cambridge University Press, Cambridge, 2004, pp.~69--90.

\bibitem[CTHS05]{T^c}
\bysame, \emph{Compactification \'{e}quivariante d'un tore (d'apr\`es {B}rylinski et {K}\"unnemann)}, Expositiones Mathematicae \textbf{23} (2005), no.~2, 161--170.

\bibitem[DSW15]{DSW12}
Ulrich Derenthal, Arne Smeets, and Dasheng Wei, \emph{Universal torsors and values of quadratic polynomials represented by norms}, Mathematische Annalen \textbf{361} (2015), no.~3, 1021--1042.

\bibitem[LL25]{LL25}
Yongqi Liang and Yufan Liu, \emph{Varieties with prescribed finite unramified {B}rauer groups and subgroups precisely obstructing the {H}asse principle}, 2025, Preprint available at \href{https://arxiv.org/abs/2504.18293}{arXiv:2504.18293}.

\bibitem[NSW13]{cnf}
J.~Neukirch, A.~Schmidt, and K.~Wingberg, \emph{Cohomology of number fields}, Grundlehren der mathematischen Wissenschaften, Springer Berlin Heidelberg, 2013.

\bibitem[VAV12]{VAV-chatletp-fold}
A.~V\'arilly-Alvarado and B.~Viray, \emph{Higher dimensional analogues of {C}h{\^a}telet surfaces.}, Bull. Lond. Math. Soc. \textbf{44} (2012), no.~1, 125--135.

\bibitem[VAV15]{VAV-model}
\bysame, \emph{Smooth compactifications of certain normic bundles}, European Journal of Mathematics \textbf{1} (2015), no.~2, 250--259.

\bibitem[Wei14]{Wei14}
Dasheng Wei, \emph{On the equation {$N_{K/k}(\Xi)=P(t)$}}, Proceedings of the London Mathematical Society \textbf{109} (2014), no.~6, 1402--1434.

\bibitem[Wei22]{wei2022+}
\bysame, \emph{The unramified {B}rauer groups of normic bundles}, 2022, Preprint available at \href{https://arxiv.org/abs/2211.07054}{arXiv:2211.07054}.

\end{thebibliography}

 \end{document}